\documentclass[10pt]{amsart}
\usepackage[T1]{fontenc}
\usepackage[utf8]{inputenc}
\usepackage[english]{babel}
\usepackage{csquotes} 
\usepackage[a4paper,margin=3cm]{geometry}
\allowdisplaybreaks
\usepackage[draft]{graphicx}
\usepackage[backend=biber,style=alphabetic,bibencoding=utf8,language=auto,autolang=other,giveninits=true,doi=false,isbn=false,url=false,maxnames=10,sorting=nty]{biblatex}
\usepackage{caption} 
\usepackage[hyperfootnotes=false]{hyperref} 
\hypersetup{
	colorlinks = true,
	linkcolor = {blue},
	urlcolor = {red},
	citecolor = {blue}
}
\usepackage[nameinlink,capitalise,sort]{cleveref}
\crefname{equation}{}{} 
\crefname{enumi}{}{} 

\crefname{figure}{Figure}{Figures}

\usepackage[shortlabels]{enumitem}
\newlist{thmenum}{enumerate}{1} 
\setlist[thmenum]{label=\textup{(\roman*)},
                  ref=\thetheorem-\textup{(\roman*)}}
\crefname{thmenumi}{Theorem}{Theorem}

\usepackage{amsmath}
\usepackage{amsthm}
\usepackage{amssymb}
\usepackage{esint} 
\usepackage{mathrsfs} 
\usepackage{faktor} 
\usepackage{dsfont} 
\usepackage[dvipsnames]{xcolor}
\usepackage{array}
\usepackage{hhline}
\usepackage[normalem]{ulem}
\usepackage{enumitem} 
\usepackage{comment} 
\usepackage[toc,page]{appendix} 
\usepackage{imakeidx} 
\usepackage{xparse} 
\usepackage{mathtools}
\usepackage{fancyhdr} 
\usepackage{ifthen} 
\usepackage{forloop} 
\usepackage{xstring}
\usepackage{tikz}
\usetikzlibrary{positioning, arrows.meta}
\usetikzlibrary{babel} 
\usetikzlibrary{cd} 
\usepackage{tikz-cd} 
\usepackage{emptypage} 
\usepackage{listings} 
\usepackage{mathabx} 
\usepackage{subcaption}

\theoremstyle{plain}
\newtheorem{lemma}{Lemma}[section]

\newtheorem{theorem}[lemma]{Theorem}

\newtheorem{remark}[lemma]{Remark}

\theoremstyle{definition}

\theoremstyle{remark}

\numberwithin{equation}{section}

\newcommand{\R}{\mathbb{R}}

\mathchardef\emptyset="001F
\renewcommand{\d}{\mathrm{d}}

\newcommand{\mynorm}{{\vert\kern-0.25ex\vert\kern-0.25ex\vert}}

\newcommand{\pt}{\partial_t}
\newcommand{\px}{\partial_x}

\newcommand{\eps}{\varepsilon}
\begin{document}


\title[Volterra equation approach]{A Volterra equation approach to the local limit of nonlocal traffic models}

\author[N.~De Nitti]{Nicola De Nitti}
\address[N.~De Nitti]{Università di Pisa, Dipartimento di Matematica, Largo Bruno Pontecorvo 5, 56127 Pisa, Italy.}
\email[]{nicola.denitti@unipi.it}

\author[K.~Huang]{Kuang Huang}
\address[K.~Huang]{The Chinese University of Hong Kong, Department of Mathematics, Shatin, New Territories, Hong Kong SAR, P.\,R.~China.}
\email[]{kuanghuang@cuhk.edu.hk}


\begin{abstract}
We consider a class of nonlocal conservation laws modeling traffic flow, given by  
\( \partial_t u_\varepsilon + \partial_x(V(u_\varepsilon \ast \gamma_\varepsilon)\, u_\varepsilon) = 0 \)  
with \( \gamma_\varepsilon(\cdot) \coloneqq \eps^{-1}\gamma(\cdot/\eps) \) for a suitable convex convolution kernel $\gamma$. Since the work of Colombo et al.~(\emph{Arch.~Ration.~Mech.~Anal.}, 2023), thanks to uniform \( \mathrm{L}^\infty \)- and TV-estimates, it is known that \( w_\varepsilon \coloneqq u_\varepsilon \ast \gamma_\varepsilon \) converges to the entropy solution of the local scalar conservation law  
\( \partial_t u + \partial_x(V(u)\, u) = 0 \) as $\varepsilon \searrow 0$.  
However, the convergence of \( \{u_\varepsilon\}_{\varepsilon > 0} \) itself has not been fully addressed so far. In this direction, a known result applies specifically to the case of an exponential kernel, where the identity \( \varepsilon \partial_x w_\varepsilon = w_\varepsilon - u_\varepsilon \) is fundamental. In this work, we address this gap in the literature and prove that \( \{u_\varepsilon\}_{\varepsilon > 0} \) converges to the same limit $u$ under the mild additional assumption that the initial datum belongs to $\mathrm L^1(\R)$. Our analysis exploits, through a Fourier approach, the stability properties of the more general Volterra-type equation \(\partial_xw_\eps=\gamma'_\eps\ast u_\eps-\gamma_\eps(0)u_\eps\), thereby deducing the convergence of $u_\eps$ from that of $w_\eps$.
\end{abstract}

\maketitle

\section{Introduction}
\label{sec:intro}

We study the following Cauchy problem for a class of nonlocal conservation laws modeling traffic flows:
\begin{align}\label{eq:cle}
\begin{cases}
\partial_t u_\eps(t,x)  + \partial_x\big(V\big((u_\varepsilon \ast \gamma_\varepsilon)(t,x)\big)\,u_\eps(t,x)\big)   	= 0,	& t >0, \ x \in \R, \\
u_\eps(0,x) = u_0(x),	&  x\in\R,
\end{cases}
\end{align}
with a parameter $\varepsilon >0$, an initial datum
\begin{align}\label{ass:u0}
    u_0 \in \mathrm  L^\infty(\R), \qquad 0 \le u_0 \le 1, \qquad \mathrm{TV}(u_0) < +\infty,
\end{align}
with $\mathrm{TV}(u_0)$ denoting the total variation of $u_0$ on $\R$, a \textit{velocity function} $V:\R \to \R_+$ satisfying 
\begin{align}\label{ass:V}
    V \in  \mathrm{Lip}([0,1]), \qquad V'(\xi) \le 0 \ \text{ for $\xi \in [0,1]$}, 
\end{align}
and a \emph{nonlocal impact} 
\[
w_\eps \coloneqq u_\eps \ast \gamma_\eps,
\]
affecting the velocity, given by a weighted average (with respect to the space variable) of the downstream density with kernel $\gamma_\eps(\cdot) \coloneqq \eps^{-1} \gamma(\cdot/\eps)$, such that 
\begin{align}\label{ass:gamma}
\begin{aligned}
&\gamma \in \mathrm{BV}(\R), \quad \operatorname{supp}\gamma \subset ]-\infty,0], \quad \int_{-\infty}^0 \gamma(z)\,\d z=1, \quad \quad \int_{-\infty}^0 |z|\gamma(z)\,\d z<+\infty, \\  
&\gamma \ge 0, \quad  \gamma \text{ convex and non-decreasing in $]-\infty,0]$}.
\end{aligned}
\end{align} 
The nonlocal impact $w_\eps$ satisfies the following evolution equation (in the strong sense): 
\begin{align} \label{eq:W}
  \pt w_\eps + \partial_x (V(w_\eps)w_\eps) =  \px \big(V(w_\eps)w_\eps - (V(w_\eps) u_\eps) \ast \gamma_\eps\big),
\end{align}
i.\,e., a conservation law with a \emph{local flux} and a \emph{nonlocal source} (in divergence form) that acts as a regularization term. 

Under assumptions \crefrange{ass:u0}{ass:gamma}, problem \eqref{eq:cle} admits a unique weak solution $u_\eps\in\mathrm L^\infty(\R_+\times\R)\cap\mathrm{C}(\R_+,\mathrm{L}^1_{\mathrm{loc}}(\R))$ and it satisfies $0 \leq u_\eps\leq 1$ (see, e.\,g.,  \cite{zbMATH06756308,MR4110434,zbMATH07615111}). 

The singular limit of \eqref{eq:cle} as $\varepsilon \searrow 0$ has attracted considerable attention, while most positive results in this direction addressed convergence of the nonlocal impact $w_\eps$ rather than the solution $u_\eps$ as the former enjoys better regularity and stability properties. For example, under assumptions \crefrange{ass:u0}{ass:gamma}, \cite{MR4553943} established the uniform total variation bound
\[ \mathrm{TV}(w_\eps(t,\cdot)) \leq \mathrm{TV}(u_0) \quad \text{for all } t \ge 0 \text{ and } \eps>0, \]
and the convergence of $w_\eps$ to the (unique) \emph{entropy solution} $u$ of the (local) scalar conservation law
\begin{align}\label{eq:cl}
\begin{cases}
    \partial_t u(t,x) + \partial_x \left(V(u(t,x)) u(t,x) \right) = 0, & t >0, \ x \in \R, \\ 
    u(0,x) = u_0(x), & x \in \R,
    \end{cases}
\end{align}
strongly in $\mathrm{L}^1_{\mathrm{loc}}(\R_+\times\mathbb{R})$ as $\varepsilon \searrow 0$. Under the same assumptions, in \cite{DeNittiKuang2025}, an asymptotically compatible Godunov-type numerical scheme was built. 

More recently, the assumptions on the total variation of the initial datum and on the convexity of the kernel have been relaxed in \cite{NCC2025} by relying on the theory of compensated compactness.

However, the convergence of $u_\eps$ in the singular limit has only been addressed in special cases, e.\,g.,  with the exponential kernel $\gamma(\cdot) \coloneqq \mathds{1}_{]-\infty,0]}(\cdot)\exp(\cdot)$ (see \cite[Theorems 3.2 and 4.2]{MR4651679}),
leveraging the identity  
\begin{align}\label{eq:exp-identity}
\varepsilon\partial_x w_\varepsilon = w_\varepsilon - u_\varepsilon.
\end{align}
In this case, the combination of \cref{eq:exp-identity} and the \( \mathrm{TV} \)-estimate on $w_{\varepsilon}$ further allows one to deduce the convergence of the family \( \{u_\varepsilon\}_{\varepsilon > 0} \) to the same limit as \( \{w_\varepsilon\}_{\varepsilon > 0} \) (see \cite[Corollary 4.1]{MR4651679}) as $\eps\searrow0$. 

Under more restrictive assumptions on the initial datum, uniform TV-bounds on $\{u_\varepsilon\}_{\varepsilon >0}$ can also be established directly: namely, when the initial datum $u_0$ is monotone, as in \cite{MR3944408}; when $u_0$ is bounded away from zero and the kernel is exponential (see \cite{MR4110434,MR4283539}); or when $u_0$ is bounded away from zero, one-sided Lipschitz continuous, and the kernel grows super-exponentially (see \cite{MR4300935}). In more general cases, uniform TV-bounds for $\{u_\varepsilon\}_{\varepsilon >0}$ are challenging (or false; see \cite{MR4300935} for some counterexamples).

In this work, we address this gap in the literature. Under the additional assumption 
\begin{align}\label{eq:L2-datum} 
u_0 \in \mathrm{L}^1(\R) ,
\end{align}
we deduce the convergence of $u_\eps(t,\cdot)$ to $u(t,\cdot)$ in $\mathrm{L}^2(\R)$ for a.\,e. $t>0$  from the convergence of $w_\eps(t,\cdot)$ to $u(t,\cdot)$.

\begin{theorem}[Convergence of $\{u_\varepsilon\}_{\varepsilon >0}$]\label{th:main}
   Let us assume that the initial datum $u_0$ satisfies \cref{ass:u0} and \cref{eq:L2-datum}, the velocity function $V$ satisfies \cref{ass:V}, and the kernel $\gamma$ satisfies \cref{ass:gamma}. Let $u_\varepsilon$ be the (unique) weak solution of \cref{eq:cle} and $w_\varepsilon \coloneqq u_\eps \ast \gamma_\eps$  the corresponding nonlocal impact that satisfies \cref{eq:W}.
   Then, for a.\,e. $t>0$, $u_\eps(t,\cdot)$ converges to $u(t,\cdot)$ in $\mathrm{L}^2(\R)$ as $\varepsilon \searrow 0$, where $u$ is the entropy solution of \cref{eq:cl}.
\end{theorem}

As a by-product of the proof of \cref{th:main}, we also obtain an explicit $\mathrm{L}^2$-convergence rate (thereby yielding convergence in $\mathrm{L}^1_{\mathrm{loc}}(\R_+\times\R)$ as well); see \cref{rmk:rate}.

\section{Proof of the main result}

The main idea for the proof of \cref{th:main} is to rely on the Volterra-type equation (of the second kind):
\[
  \partial_xw_\eps(t,x)=(\gamma'_\eps\ast u_\eps)(t,x)-\gamma_\eps(0)u_\eps(t,x) = -\gamma_\eps(0)u_\eps(t,x) + \int_x^\infty u_\eps(t,y)\gamma'_\eps(x-y)\,\d y,
\]
which generalizes the identity \cref{eq:exp-identity}. We exploit stability properties of this Volterra-type equation through a Fourier approach (a strategy previously used in \cite{du2018stability,huang2022stability}).

\begin{proof}[Proof of \cref{th:main}]
Owing to \cite[Theorem 1.3]{MR4553943}, under the assumptions \crefrange{ass:u0}{ass:gamma}, we have 
\[w_\varepsilon \to u \quad \text{ in $\mathrm{L}^1_{\mathrm{loc}}(\R_+ \times \R)$ as $\varepsilon \searrow 0$},\]
where $\{w_\varepsilon\}_{\varepsilon >0}$ is the family of solutions to \cref{eq:W} (with initial datum $u_0 \ast \gamma_\varepsilon$) and $u$ is the unique entropy solution of \cref{eq:cl}. We aim to prove that, under the additional assumption \cref{eq:L2-datum}, $\{u_\varepsilon\}_{\varepsilon >0}$ converges to the same limit in $\mathrm L^2([0,T] \times \R)$ for any fixed $T>0$.

\uline{Step 1.} \emph{Energy estimate for $u_\varepsilon$.}
Fom  the identity $w_\eps=\gamma_\eps \ast u_\eps$, we deduce that  
\begin{align}\label{eq:dx_w}
    \partial_xw_\eps=\gamma'_\eps\ast u_\eps-\gamma_\eps(0)u_\eps,
\end{align}
and thus, $\|\partial_xw_\eps\|_{\mathrm{L}^\infty} \leq \gamma_\eps(0)$.
Then, the estimate
\begin{align*}
    \frac{\d}{\d t} \int_\R |u_\eps(t,x)|^2\,\d x &= -\int_\R |u_\eps(t,x)|^2\,V'(w_\eps(t,x))\,\partial_xw_\eps(t,x)\,\d x \\ &\leq \|V'\|_{\mathrm{L}^\infty}\, \gamma_\eps(0) \int_\R |u_\eps(t,x)|^2\,\d x,
\end{align*}
gives that $u_\eps(t,\cdot)\in\mathrm{L}^2(\R)$ for any $t\geq0$, and consequently $w_\eps(t,\cdot)\in\mathrm{L}^2(\R)$ for any $t\geq0$. This allows us to consider the singular limit problem in the Fourier domain. 

\uline{Step 2.} \emph{A Fourier approach.} Applying the Fourier transform (in $\mathrm L^2(\R)$) on $w_\eps=\gamma_\eps\ast u_\eps$, we get $\hat{w}_\eps(t,\cdot)=\hat{\gamma}_\eps\cdot\hat{u}_\eps(t,\cdot)$ for all $t\geq0$, where
\begin{align*}
    &\hat{\gamma}_\eps(\xi) = \int_{-\infty}^0 e^{-\mathrm{i}\eps\xi z} \gamma(z)\,\d z = a_\eps(\xi)+\mathrm{i} \, b_\eps(\xi), \\
    &\textnormal{with} \quad  a_\eps(\xi) \coloneqq \int_{-\infty}^0 \cos(\eps\xi z) \gamma(z) \,\d z \quad \textnormal{and} \quad b_\eps(\xi) \coloneqq -\int_{-\infty}^0 \sin(\eps\xi z) \gamma(z) \,\d z.
\end{align*}
The assumption \cref{ass:gamma} yields that
\begin{align*}
    b_\eps(\xi) = \frac{1}{\eps\xi} \int_{-\infty}^0 (1-\cos(\eps\xi z)) \gamma'(z)\,\d z > 0 \quad \text{for all } \eps>0,\xi\neq0.
\end{align*}
When $\xi=0$, $b_\eps(\xi)=0$ but $a_\eps(\xi)=1$. Thus, we have $\hat{\gamma}_\eps(\xi)\neq0$ for all $\eps>0$ and $\xi\in\R$. Moreover, we have
\[ \hat{\gamma}_\eps(\xi) = 1 + \mathrm{i}\eps\xi \int_{-\infty}^0  \frac{1-e^{-\mathrm{i}\eps\xi z}}{\mathrm{i}\eps\xi z}\cdot |z|\gamma(z)\,\d z \sim1 +\mathrm{i}\eps\xi\int_{-\infty}^0 |z|\gamma(z)\,\d z \to1 \quad \text{ as }\eps\xi\to0. \]
Therefore, we deduce the following pointwise convergence from the $\mathrm{L}^2$-convergence of $\hat{w}_\eps$ to $\hat{u}$:
\begin{align*}
    \hat{u}_\eps(t,\xi)=\frac{\hat{w}_\eps(t,\xi)}{\hat{\gamma}_\eps(\xi)} \to \hat{u}(t,\xi) \quad \textnormal{as }\eps\searrow0\textnormal{ for a.\,e. }t\geq0,\ \xi\in\R.
\end{align*}

We need, however, to estimate the $\mathrm{L}^2$-distance $\|u_\eps(t,\cdot)-u(t,\cdot)\|_{\mathrm{L}^2(\R)}$. To this end, using Plancherel's theorem and splitting into high and low frequencies, we compute as follows (for a given $\delta >0$ to be chosen later):  
\begin{align*}
    \|u_\eps(t,\cdot)-u(t,\cdot)\|_{\mathrm{L}^2(\R)}^2 &= \int_\R \left|\hat{u}_\eps(t,\xi) - \hat{u}(t,\xi)\right|^2 \,\d \xi \\
    &= \int_{|\eps\xi|<\delta} \left|\frac{1}{\hat{\gamma}_\eps(\xi)} \hat{w}_\eps(t,\xi) - \hat{u}(t,\xi)\right|^2 \,\d \xi + \int_{|\eps\xi|\geq\delta} \left|\hat{u}_\eps(t,\xi) - \hat{u}(t,\xi)\right|^2 \,\d \xi \\
    &\leq \int_{|\eps\xi|<\delta} \frac{1}{|\hat{\gamma}_\eps(\xi)|^2} \left| \hat{w}_\eps(t,\xi) - \hat{u}(t,\xi) \right|^2 \,\d \xi + \int_{|\eps\xi|<\delta} \left|\frac{1}{\hat{\gamma}_\eps(\xi)} - 1\right|^2 \left|\hat{u}(t,\xi)\right|^2 \,\d \xi \\
    &\qquad + 2 \int_{|\eps\xi|\geq\delta} \left|\hat{u}_\eps(t,\xi)\right|^2 \,\d \xi + 2 \int_{|\eps\xi|\geq\delta} \left|\hat{u}(t,\xi)\right|^2 \,\d \xi \\
    &\eqqcolon \mathcal{I}_1+\mathcal{I}_2+\mathcal{I}_3+\mathcal{I}_4.
\end{align*}

\uline{Step 3.} \emph{Estimating the low frequencies  (terms $\mathcal I_1$ and $\mathcal I_2$).}  Using the fact that $\hat{\gamma}_\eps(\xi)\sim1 +\mathrm{i}\eps\xi\int_{-\infty}^0 |z|\gamma(z)\,\d z$ as $\eps\xi\to0$, by taking $\delta$ small enough, we have $|\hat{\gamma}_\eps(\xi)|\geq\frac12$ when $|\eps\xi|<\delta$ and 
\begin{align*}
\sup_{|\eps\xi|<\delta} \left|\frac{1}{\hat{\gamma}_\eps(\xi)} - 1\right| \leq C_0\delta,
\end{align*}
where the constant $C_0>0$ depends only on $\gamma$. Thus, using Plancherel's theorem and the uniform bounds $0\leq w_\eps,u\leq1$, we deduce that
\begin{align*}
    \mathcal{I}_1 \leq& 4 \int_{|\eps\xi|<\delta} \left| \hat{w}_\eps(t,\xi) - \hat{u}(t,\xi) \right|^2 \,\d \xi \leq 4\|w_\eps(t,\cdot)-u(t,\cdot)\|_{\mathrm{L}^2(\R)}^2 \leq 4 \|w_\eps(t,\cdot)-u(t,\cdot)\|_{\mathrm{L}^1(\R)}, \\
    \mathcal{I}_2 \leq& C_0^2\delta^2 \int_{|\eps\xi|<\delta} \left|\hat{u}(t,\xi)\right|^2 \,\d \xi \leq C_0^2\delta^2 \|u(t,\cdot)\|_{\mathrm{L}^2(\R)}^2 \leq C_0^2\delta^2 \|u(t,\cdot)\|_{\mathrm{L}^1(\R)} = C_0^2\delta^2 \|u_0\|_{\mathrm{L}^1(\R)}.
\end{align*}

\uline{Step 4.} \emph{Estimating the high frequencies (terms $\mathcal I_3$ and $\mathcal I_4$).} To estimate the terms 
$\mathcal I_3$ and $\mathcal I_4$, we derive further estimates on the family $\{\hat{u}_\eps\}_{\eps>0}$. Multiplying \cref{eq:dx_w} by $u_\eps$ and integrating it over $x\in\R$, we obtain
\begin{align*}
    \int_\R u_\eps(t,x)\partial_xw_\eps(t,x)\,\d x &= \int_{-\infty}^0 \left( \int_\R u_\eps(t,x)u_\eps(t,x-s) - u_\eps(t,x)^2 \,\d x \right) \gamma'_\eps(s) \,\d s \\
    &= -\frac12 \int_{-\infty}^0 \left( \int_\R \left( u_\eps(t,x-s) - u_\eps(t,x) \right)^2 \,\d x \right) \gamma'_\eps(s) \,\d s,
\end{align*}
for all $t\geq0$. Letting
\[ \rho_\eps(t,s)\coloneqq \int_\R \left( u_\eps(t,x-s) - u_\eps(t,x) \right)^2 \,\d x \quad \text{for all $t\geq0$, $s\in\R$, and $\eps>0$}, \]
we deduce that
\begin{align}\label{eq:rhoeps_estm}
    0 \leq \int_{-\infty}^0 \rho_\eps(t,s) \gamma'_\eps(s) \,\d s \leq 2 \left|\int_\R u_\eps(t,x)\partial_xw_\eps(t,x)\,\d x\right| \leq 2\,\mathrm{TV}(w_\eps(t,\cdot)) \leq 2\,\mathrm{TV}(u_0),
\end{align}
for all $t\geq0$ and $\eps>0$.

On the other hand, in the Fourier domain, we have
\begin{align*}
    \rho_\eps(t,s) = \int_\R \left| e^{-\mathrm{i} s\xi} -1\right|^2 \left|\hat{u}_\eps(t,\xi)\right|^2 \,\d \xi = \int_\R 2(1-\cos(s\xi)) \left|\hat{u}_\eps(t,\xi)\right|^2 \,\d \xi,
\end{align*}
which yields
\begin{align*}
    \int_{-\infty}^0 \rho_\eps(t,s) \gamma'_\eps(s) \,\d s =& \int_\R \left( \int_{-\infty}^0 2(1-\cos(s\xi)) \gamma'_\eps(s) \,\d s \right) \left|\hat{u}_\eps(t,\xi)\right|^2 \,\d \xi \\
    =& \frac1\eps \int_\R \left( \int_{-\infty}^0 2(1-\cos(\eps s\xi)) \gamma'(s) \,\d s \right) \left|\hat{u}_\eps(t,\xi)\right|^2 \,\d \xi.
\end{align*}
We define
\begin{align*}
    h(z) \coloneqq \int_{-\infty}^0 2(1-\cos(z s)) \gamma'(s) \,\d s \quad \text{for } z\in\R.
\end{align*}
Noting that $\gamma'\geq0$ and $\int_{-\infty}^0\gamma'(z)\,\d z=\gamma(0)>0$, we have $h(z)>0$ for all $z\neq0$. 
When $z\to0$, the mean value theorem yields
\[ 2(1-\cos(z s)) = \cos(\zeta(z s))(z s)^2 \quad \text{for } \zeta(z s) \text{ between } 0 \text{ and } z s, \]
thus
\[ \left| \frac{2(1-\cos(zs))}{z^2} \right| \leq s^2 \quad \text{and} \quad \lim_{z\to0} \frac{2(1-\cos(zs))}{z^2} = s^2 \quad \text{for all }s<0. \]
Therefore, by Lebesgue's dominated convergence theorem, we have
\[ \lim_{z\to0} \frac{h(z)}{z^2} = \int_{-\infty}^0 s^2\gamma'(s)\,\d s = 2\int_{-\infty}^0 |s|\gamma(s)\,\d s \in ]0,+\infty[. \]
Moreover, Riemann--Lebesgue's lemma yields that
\[ \lim_{|z|\to\infty} h(z) = 2\int_{-\infty}^0 \gamma'(s) \,\d s=2\,\gamma(0)\in ]0,+\infty[. \]
By continuity of $h$ and the two limits established above, there exist $\delta_\gamma,\eta_\gamma>0$, depending only on $\gamma$, such that
\[ h(z)\geq \eta_\gamma z^2\geq\eta_\gamma \delta^2 \, \text{ for }\, |z|\geq\delta \quad \text{whenever } 0<\delta<\delta_\gamma. \]
Then we take $z=\eps\xi$ and obtain from \cref{eq:rhoeps_estm} that
\[ 2\,\mathrm{TV}(u_0) \geq \frac1\eps\int_\R h(\eps\xi)\left|\hat{u}_\eps(t,\xi)\right|^2 \,\d \xi \geq \frac{\eta_\gamma \delta^2}{\eps} \int_{|\eps\xi|\geq\delta} \left|\hat{u}_\eps(t,\xi)\right|^2 \,\d \xi, \]
and thus
\begin{align}\label{eq:unifm_int}
    \int_{|\eps\xi|\geq\delta} \left|\hat{u}_\eps(t,\xi)\right|^2 \,\d \xi \leq \frac{2\eps}{\eta_\gamma \delta^2}\mathrm{TV}(u_0) \eqqcolon C_1 \,\mathrm{TV}(u_0)\frac{\eps}{\delta^2},
\end{align}
for any $\delta\in]0,\delta_\gamma[$. This estimates $\mathcal I_3$. 

For the term $\mathcal{I}_4$, we use the TV-estimate $\mathrm{TV}(u(t,\cdot))\leq\mathrm{TV}(u_0)$ for all $t\geq0$ to deduce that
\begin{align*}
    \int_{|\eps\xi|\geq\delta} \left|\hat{u}(t,\xi)\right|^2 \,\d \xi \leq \int_{|\eps\xi|\geq\delta} \left(\frac{\mathrm{TV}(u(t,\cdot))}{|\xi|}\right)^2\,\d\xi = 2 (\mathrm{TV}(u_0))^2 \,\frac{\eps}{\delta}.
\end{align*}

\uline{Step 5.} \emph{Conclusion of the proof.} Combining the estimates in Steps 3--4, we obtain
\begin{equation}\label{eq:convergence}
\begin{aligned}
    \|u_\eps(t,\cdot)-u(t,\cdot)\|_{\mathrm{L}^2(\R)}^2 \leq& 4 \,\|w_\eps(t,\cdot)-u(t,\cdot)\|_{\mathrm{L}^1(\R)} + C_0^2\delta^2 \,\|u_0\|_{\mathrm{L}^1(\R)} \\
    &\qquad + C_1\,\mathrm{TV}(u_0)\frac{\eps}{\delta^2}+ 2 (\mathrm{TV}(u_0))^2 \,\frac{\eps}{\delta}.
\end{aligned}
\end{equation}
Letting $\eps\searrow0$ and $\delta=\eps^{1/4}$, we have $\|w_\eps(t,\cdot)-u(t,\cdot)\|_{\mathrm{L}^1(\R)}\to0$ and $\delta^2,\,\eps/\delta^2,\,\eps/\delta\searrow0$, giving $\|u_\eps(t,\cdot)-u(t,\cdot)\|_{\mathrm{L}^2(\R)}\to0$.
\end{proof}

\begin{remark}[Convergence rate]\label{rmk:rate}
In \cite{MR4553943}, the convergence rate
\begin{align}\label{eq:convr-w-e} \|w_\eps(t,\cdot)-u(t,\cdot)\|_{\mathrm{L}^1(\R)} \leq C(\gamma,V)\,\mathrm{TV}(u_0) \left(\eps+\sqrt{\eps t}\right) 
\end{align}
was established. 
By taking $\delta=\eps^{1/4}$ in \cref{eq:convergence} and using \cref{eq:convr-w-e}, we deduce
\[ \|u_\eps(t,\cdot)-u(t,\cdot)\|_{\mathrm{L}^2(\R)} \leq C(\gamma,V)\,\|u_0\|_{\mathrm{BV}(\R)}^{\frac12} \left(\eps(1+t)\right)^{\frac14}, \]
provided $0<\eps<\eps_0\coloneqq\left(1+\mathrm{TV}(u_0)\right)^{-2}$, where $\|u_0\|_{\mathrm{BV}(\R)}\coloneqq\|u_0\|_{\mathrm{L}^1(\R)}+\mathrm{TV}(u_0)$ and the constant $C(\gamma,V)>0$ depends only on the kernel $\gamma$ and the velocity function $V$.
We leave the question of the sharpness of this rate (as well as the search for optimal exponents) for future investigation.
\end{remark}

\section*{Acknowledgments}

N.~De Nitti is a member of the Gruppo Nazionale per l’Analisi Matematica, la Probabilità e le loro Applicazioni (GNAMPA) of the Istituto Nazionale di Alta Matematica (INdAM). He thanks The Chinese University of Hong Kong, where part of this work was carried out, for their kind hospitality. 

K.~Huang was supported by a Direct Grant of Research (2024/25) from The Chinese University of Hong Kong.

We thank G.~M.~Coclite and B.~Jin for their helpful and constructive comments.

\vspace{0.5cm}

\printbibliography

\vfill 

\end{document}